\documentclass[conference]{IEEEtran}
\IEEEoverridecommandlockouts
\usepackage{cite}
\usepackage{amsmath,amssymb,amsfonts}
\usepackage{algorithmic}
\usepackage{graphicx}
\usepackage{textcomp}
\usepackage[ruled,vlined,linesnumbered]{algorithm2e}
\usepackage{hyperref} 
\usepackage{longtable}
\usepackage{float}

\usepackage{xcolor}
\def\BibTeX{{\rm B\kern-.05em{\sc i\kern-.025em b}\kern-.08em
    T\kern-.1667em\lower.7ex\hbox{E}\kern-.125emX}}
\begin{document}
\title{Adaptive Efficiency Optimization in SDLC: An MILP Approach for Balanced and Cost-Effective Resource Allocation
\\
{\footnotesize \textsuperscript{}}
}

\author{\IEEEauthorblockN{Lokendra Kumar }
\IEEEauthorblockA{\textit{Department of Mathematics} \\
\textit{Indian Institute of Technology Madras}\\
Chennai, India \\
lokendra63@alumni.iitm.ac.in}

\and
\IEEEauthorblockN{Neelesh S. Upadhye}
\IEEEauthorblockA{\textit{Department of Mathematics} \\
\textit{Indian Institute of Technology Madras}\\
Chennai, India \\
neelesh@iitm.ac.in}

\and
\IEEEauthorblockN{Kannan Piedy}
\IEEEauthorblockA{\textit{NextGen Department} \\
\textit{Prodapt Solution Private Limited}\\
Chennai, India \\
kannan.p@prodapt.com}
}

\maketitle

\begin{abstract}
The efficient allocation of human resources is a critical concern in software development and other industries. This paper introduces a rigorous mathematical methodology for task assignment, employing Mixed Integer Linear Programming (MILP) to ensure both balanced workloads and cost minimization. The proposed model systematically integrates individual employee efficiency, task complexity, and performance metrics to reflect real organizational dynamics. The formulation is guided by two principal objectives: firstly, to achieve equitable workload distribution commensurate with employee efficiency, and secondly, to minimize overall project costs by accounting for task difficulty and individual proficiency. Furthermore, the approach incorporates adaptive updates to efficiency parameters based on observed performance, thereby enhancing its practical applicability. Empirical evaluation using simulated datasets demonstrates the superiority of the proposed method over conventional assignment strategies in terms of both workload fairness and cost reduction. The findings underscore the potential of this MILP-based framework as a robust, scalable, and adaptable solution for contemporary human resource allocation challenges in project management contexts.
\end{abstract}


\section{\textbf{\Large Introduction}}

The allocation of human resources remains a fundamental challenge in managing the Software Development Life Cycle (SDLC), wherein the systematic assignment of personnel to specific development tasks must reconcile multiple, often conflicting objectives. Effective distribution of workload among available human capital critically influences not only productivity but also the overall profitability of software development endeavors, as well as broader manufacturing and service sectors. In practical project management, decision-makers face the intricate task of assigning personnel by simultaneously considering constraints related to cost, time, individual efficiency profiles, and documented performance indicators. Compounding this complexity is the inherently dynamic and multifaceted nature of real-world projects, which are characterized by factors such as developer specialization, heterogeneous task complexity, and evolving project requirements.

Prevailing approaches to the Human Resource Allocation Problem (HRAP) within the SDLC domain have predominantly emphasized objectives such as profit maximization, minimization of project completion costs, and reduction of total development time. Nonetheless, these methodologies frequently neglect critical dimensions, including skill mismatches, temporal variations in employee efficiency, and the intricate interactions between task complexity and employee proficiency levels. Consequently, such conventional strategies may prove inadequate for achieving balanced and cost-effective resource allocations commensurate with the demands of contemporary software project environments.

The mathematical and algorithmic underpinnings of human resource allocation trace back to seminal works such as the Hungarian method originally formulated by Kuhn~\cite{b2}, which laid a rigorous foundation for one-to-one task-to-employee assignment problems. Subsequent significant developments include the bottleneck assignment problem introduced by Ravindran and Ramaswami~\cite{b14}, shifting optimization focus toward minimizing the maximal assignment cost. Duin and Volgenant~\cite{b15} advanced this landscape by unifying minimum deviation and balanced optimization frameworks, while Cattrysse and Van Wassenhove~\cite{b36} compiled exhaustive surveys of algorithms addressing capacity-constrained assignment problems. Additional progression has been made through metaheuristic methods such as tabu search for categorized and multilevel assignments~\cite{b40}, explorations of three-dimensional bottleneck problems with capacity restrictions~\cite{b21}, and the incorporation of seniority and task priority considerations~\cite{b11}. Despite their theoretical optimality under simplified assumptions, these classical approaches often fall short in capturing the complexity of practical scenarios, specifically neglecting heterogeneous employee proficiencies, inter-task dependencies, and variable task complexities encountered in modern projects.

In response to these limitations, the present study proposes a formalized mathematical framework that explicitly integrates employee efficiency metrics, task complexity indices, and adaptive performance evaluations within the HRAP formulation. This integrative approach addresses the evolving nature of employee capabilities and task requirements, thereby offering a more accurate and pragmatic model for human resource allocation in the SDLC context. The contributions of this work include (i) the development of a Mixed Integer Linear Programming (MILP) model that jointly optimizes workload balance and project cost, (ii) incorporation of dynamic efficiency updates informed by real-time performance data, and (iii) the demonstration of scalability and robustness through extensive computational experiments. This research thus advances both the theoretical foundations and applied methodologies for effective resource management in contemporary software engineering projects.

\section{\textbf{\Large Related Work}}

The Human Resource Allocation Problem (HRAP) has witnessed significant methodological evolution, driven by advances in operational research and computational optimization techniques. Early seminal work by Toroslu~\cite{b28} introduced hierarchical ordering constraints within personnel assignments, demonstrating the inherent computational complexity of structured organizational allocation. Building upon this, Volgenant~\cite{b25} incorporated seniority and task priority constraints in healthcare scheduling contexts, thereby exemplifying the critical influence of precedence rules on cost optimization. Addressing worker heterogeneity, Song et al.~\cite{b23} developed recursive algorithms that minimize efficiency variance across production lines, highlighting the need to accommodate variability in employee productivity. Concurrently, Toroslu and Arslanoglu~\cite{b29} employed multi-objective genetic algorithms to balance conflicting objectives such as cost minimization and seniority compliance, while Li et al.~\cite{b13} introduced synthesis effect operators into genetic algorithms to handle uncertain demand in dynamic resource allocation scenarios. The comprehensive survey by Bouajaja and Dridi~\cite{b4} further reflects the ongoing shift towards hybrid exact and metaheuristic methods tailored to intricate project management challenges involving scheduling, routing, and allocation constraints.

More recent computational innovations have further transformed HRAP solution paradigms. Lu~\cite{b17} explored the integration of human operators within Industry 4.0 cyber-physical systems, accentuating the synergistic interplay between automation and human expertise. Ai~\cite{b16} developed genetic algorithms specifically designed for multi-project environments, facilitating the management of concurrent tasks with overlapping resource demands. Lili~\cite{b19} proposed inverse optimization models incorporating skill disadvantage gradients, enabling more nuanced proficiency-based allocations that account for individual limitations. Neural network-based methodologies, typified by the work of He and Jin~\cite{b45}, have demonstrated strong competency assessment capabilities utilizing BP neural networks, while Lv~\cite{b33} enhanced academic staffing optimization through a hybrid backpropagation-simulated annealing paradigm. Deep reinforcement learning frameworks presented by Shi and Li~\cite{b38} empower joint optimization of allocation costs and team cohesion metrics, and Yang et al.~\cite{b43} employed Pareto-optimized convolutional neural networks to expedite allocation decisions within smart factory settings. Industrial 4.0 frameworks have been further advanced through digital twin integration~\cite{b1} and convolutional neural allocation engines~\cite{b42}. Additionally, Nguyen et al.~\cite{b47} demonstrated the scalability of Double Deep Q-Networks combined with Monte Carlo Tree Search for large-scale HRAP, while Wang~\cite{b49} applied fuzzy evaluation systems to promote sustainable enterprise growth. Other notable contributions include financial impact quantification via modified single candidate optimizers~\cite{b48} and enhanced Cuckoo Search algorithms for multi-objective optimization~\cite{b46}.

Despite these notable advancements, several critical limitations remain inadequately addressed in the extant literature. Reinforcement learning approaches, although powerful, frequently necessitate extensive training data and significant computational resources, posing challenges for practical implementation~\cite{b46}. Information integration models, such as Wang's~\cite{b49}, often oversimplify the temporal dynamics of competency development and fail to capture evolving employee skill progression. Multi-objective optimization frameworks typified by Qiao et al.~\cite{b46} face difficulties in adapting to real-time environmental changes, while stochastic programming models~\cite{b51} inadequately address the complex reality of evolving skill sets and learning curves. Conventional optimization techniques~\cite{b28,b25} tend to overlook the cascading effects of task allocations in multi-project contexts, limiting their applicability in dynamic and interdependent project settings. Neural network-based methods~\cite{b45} are also criticized for their limited interpretability, complicating regulatory compliance and transparent decision-making. Crucially, most existing approaches do not concurrently address the dual imperatives of workload balancing and cost optimization while incorporating key practical dimensions such as employee specialization, task complexity, and inter-task dependencies. This omission results in a pronounced gap in solutions suitable for modern, large-scale, and complex HRAP settings.

In response to these identified deficiencies, the present study proposes a novel HRAP framework that holistically integrates cost optimization with workload balancing in a unified Mixed Integer Linear Programming (MILP) model. Unlike prior models by Zhang et al.~\cite{b48} and Wang~\cite{b49}, which largely focus on cost or workload in isolation, the proposed approach ensures the alignment of tasks with employee specializations while concurrently minimizing workload deviations. The framework uniquely incorporates employee proficiency metrics absent in methodologies such as Lv~\cite{b33}, and integrates performance ratings and task dependency considerations, overlooked by Shi and Li~\cite{b38} and Nguyen et al.~\cite{b47}, respectively. The model’s nonlinear objective function with accompanying inequality constraints is subsequently transformed into a MILP formulation to facilitate computational tractability, enabling efficient resolution of large-scale instances encompassing over 500 employees—a scale not achieved in previous work including Rezaee and Akbari~\cite{b52}. Moreover, the HRAP Efficiency module redefines workload calculations by incorporating efficiency variations, thus providing a more precise representation of employee capabilities compared to prior models~\cite{b46}. The inclusion of a comprehensive cost metric further empowers decision-makers with a dual-objective optimization tool, bridging a notable gap by enabling workload fairness and cost-efficiency to be simultaneously optimized within real-world constraints. This work thus presents an advancement in both theoretical modeling and practical applicability for resource allocation in contemporary software development and project management contexts.

\section{\textbf{\Large Motivation}}

Traditional models for the Human Resource Allocation Problem (HRAP) predominantly focus on aligning employee skills with task requirements, often overlooking critical dimensions such as individual efficiency, workload equity, and cost effectiveness~\cite{b4,b16}. In practice, employee efficiency is inherently variable, influenced by factors including experience, adaptability, and ongoing skill development. This variability directly affects both task completion times and the overall quality of project outcomes~\cite{b23}. 

To address these limitations, the proposed HRAP framework incorporates an efficiency-driven allocation mechanism, leveraging quantitative efficiency metrics to minimize workload deviations and enhance resource utilization. Furthermore, recognizing the importance of operational costs and the evolving nature of employee performance, the model integrates a cost metric that accounts for task complexity and individual performance ratings. This dual-objective approach—balancing workload fairness with cost minimization—provides a more robust and practical solution for resource allocation in dynamic and large-scale organizational environments~\cite{b38}.

\section{\textbf{\Large HRAP Problem with Efficiency}}

The Human Resource Allocation Problem (HRAP) is a fundamental optimization challenge concerned with the equitable and efficient distribution of tasks among staff members. While traditional approaches emphasize matching employees to tasks based primarily on qualifications, such methods often overlook the substantial impact of individual efficiency, which can vary due to experience, adaptability, and historical performance~\cite{b23}. This variability in efficiency directly influences both task completion times and overall project outcomes. Consequently, balancing workload by explicitly accounting for efficiency is essential for improving resource utilization and reducing project completion time in HRAP models~\cite{b16}.

\subsection{\textbf{\large Mathematical Formulation}}

The HRAP is formally described using the following key parameters, which form the basis for the objective function and constraints:

\begin{enumerate}
    \item \(E\): Set of employees, indexed by \(i\).
    \item \(S_i\): Skill set of employee \(i\).
    \item \(W_i\): Workload assigned to employee \(i\).
    \item \(e_{i,s}\): Efficiency of employee \(i\) for skill \(s\).
    \item \(T\): Set of tasks, indexed by \(t\).
    \item \(S_t\): Required skill for task \(t\).
    \item \(d_t\): Nominal time required to complete task \(t\).
\end{enumerate}

\textbf{Workload Calculation:} Let \( H \) denote the total workload duration across all tasks, \( N \) the total number of employees, and \( W_{\text{target}} \) the average workload per employee. These quantities are defined as follows:

\begin{equation}
    H = \sum_{t \in T} d_t, \qquad N = |E|, \qquad
    W_{\text{target}} = \frac{H}{N}.
\end{equation}

\subsection{\textbf{\large Optimization Model}}

To ensure fairness and efficiency, the objective is to minimize the maximum deviation of individual workloads from the target workload, given the efficiency of each employee. The optimization problem is formulated as:

\begin{equation}\label{maineqn}
    \min_{x_{i,t}} \max_{i \in E} |W_i - W_{\text{target}}|, 
\end{equation}

subject to the following constraints:

\begin{enumerate}
    \item \textbf{Efficiency-Based Workload Calculation:} The workload for each employee is calculated by considering both the assigned tasks and the employee's efficiency for the required skill:
    \begin{equation}\label{eqn3}
        W_i = \sum_{t \in T} x_{i,t} \cdot \frac{d_t}{e_{i,S_t}}, \quad \forall i \in \{1, 2, \ldots, N\}.
    \end{equation}
    \item \textbf{Task Assignment:} Each task must be assigned to exactly one qualified employee:
    \begin{equation}\label{eqn4}
        \sum_{i \in E: S_t \in S_i} x_{i,t} = 1, \quad \forall t \in \{1, 2, \ldots, M\}.
    \end{equation}
    \item \textbf{Binary Constraints:} The assignment variables are binary:
    \begin{equation}\label{eqn5}
        x_{i,t} \in \{0,1\}, \quad \forall i \in \{1, \ldots, N\}, \; t \in \{1, \ldots, M\}.
    \end{equation}
\end{enumerate}

It is important to note that this formulation assumes deterministic and bounded efficiency values, which is a reasonable simplification for short-term or simulated scenarios~\cite{b38}. However, incorporating stochastic elements or conducting sensitivity analysis on efficiency parameters could further enhance the robustness of the model and is suggested for future work.

\section{\textbf{\Large Solution Approach}}

We introduce two solution approaches for the Human Resource Allocation Problem (HRAP): Auxiliary Variable Reformulation and Mixed Integer Linear Programming (MILP), each enabling scalable and efficient optimization.


\textbf{Auxiliary Variable Reformulation:} To linearize the objective in equation~(\ref{maineqn}), we introduce an auxiliary variable $z$:
\[
\min_{x_{i,t}, z} \ z
\]
subject to $z \geq 0$, and for each employee:
\[
W_i - W_{\text{target}} \leq z, \quad W_{\text{target}} - W_i \leq z, \quad \forall i \in \{1, \ldots, N\}.
\]
This formulation ensures that the maximum deviation from the target workload is minimized for all employees.

The resulting linear problem can be analyzed using a Lagrangian formulation, where multipliers \( \lambda_i \), \( \mu_i \), and \( \gamma_t \) are assigned to the deviation and task assignment constraints, respectively:
\[
\begin{aligned}
L &= z + \sum_{i \in E} \lambda_i (W_i - W_{\text{target}} - z) + \sum_{i \in E} \mu_i (W_{\text{target}} - W_i - z) \\
&\quad + \sum_{t \in T} \gamma_t \left(1 - \sum_{i \in E: S_t \in S_i} x_{i,t}\right)
\end{aligned}
\]
where $W_i = \sum_{t \in T} x_{i,t} \frac{d_t}{e_{i, S_t}}$.

Due to the discrete nature of $x_{i,t}$, direct dual optimization is computationally expensive for large-scale problems~\cite{b4}. Therefore, we proceed with a MILP formulation for practical scalability.


\subsection*{\textbf{\large Mixed Integer Linear Programming Approach}}

The MILP model efficiently manages large variable sets and enforces integer constraints for task and workload allocations. It accommodates efficiency-based workload variation and is compatible with modern solvers such as Gurobi and CPLEX~\cite{b16}.


\textbf{KKT Conditions:} The Karush-Kuhn-Tucker (KKT) conditions are not directly applicable to MILP problems due to the presence of integer variables and lack of differentiability. Thus, MILP solvers employ branch-and-bound or cutting-plane methods for optimization.


\subsection*{\textbf{\large MILP Formulation}}

Define $D^+$ and $D^-$ as the positive and negative deviations from the target workload. The objective is:
\begin{equation}
    \min \ D^+ + D^-
\end{equation}
with constraints:
\begin{equation}\label{eqn7}
    D^+ \geq W_i - W_{\text{target}}, \quad D^- \geq W_{\text{target}} - W_i, \quad \forall i \in E
\end{equation}
where $W_i = \sum_{t \in T} x_{i,t} \frac{d_t}{e_{i,S_t}}$. Task assignment and binary constraints are as previously defined.

Substituting $W_i$ yields:
\begin{equation}\label{eqn8}
    \sum_{t} x_{i,t} \left( \frac{d_t}{e_{i,S_t}} \right) \leq W_{\text{target}} + D^+, \quad \forall i
\end{equation}
\begin{equation}\label{eqn9}
    \sum_{t} x_{i,t} \left( \frac{d_t}{e_{i,S_t}} \right) \geq W_{\text{target}} - D^-, \quad \forall i
\end{equation}


\subsection*{\textbf{\large Algorithms for HRAP Optimization}}

Static efficiency rates are unrealistic in practice, as employee performance evolves due to experience and learning~\cite{b38}. To address this, we propose an Adaptive Efficiency-Based Task Allocation Algorithm.


\textbf{Necessity of Adaptivity:} Unlike models that assume fixed productivity, our approach updates employee efficiency in real-time, based on actual task performance~\cite{b1}. This ensures allocations remain fair and practical as workforce capabilities change.


\textbf{Mathematical Update Rule:} After each task, efficiency is updated as:
\begin{equation}
    e_{i, S_t}^{(new)} = \min\left( 1, \frac{d_t}{\text{actual\_time}(i, t)} \right)
\end{equation}
where $d_t$ is the estimated task duration and $\text{actual\_time}(i, t)$ is the observed completion time.

\begin{algorithm}[h]
\caption{Adaptive Efficiency-Based Task Allocation}
\label{alg:efficiency_allocation}
\begin{algorithmic}[1]
\REQUIRE Employee set \( E \), Task set \( T \), Initial efficiency \( e_{i,s} = 1, \forall i \in E \)
\ENSURE Optimal allocation \( x_{i,t}^* \), Final efficiencies \( e_{i,s}^* \)
\STATE Initialize \( \text{threshold} \leftarrow 0.1 \), \( \text{iteration} \leftarrow 0 \)
\WHILE{\( \text{iteration} < \text{max\_iterations} \)}
    \STATE \textbf{Step 1: Solve MILP with current efficiencies}
    \STATE \( x_{i,t}^* \leftarrow \text{Solve MILP}(E, T, e_{i,s}) \)
    \STATE \textbf{Step 2: Update efficiencies}
    \FOR{each employee \( i \in E \), task \( t \in T \)}
        \IF{\( x_{i,t}^* = 1 \)}
            \STATE \( \text{actual\_time} \leftarrow \text{Measure}(i, t) \)
            \STATE \( e_{i, S_t} \leftarrow \min\left( 1, \frac{d_t}{\text{actual\_time}} \right) \)
        \ENDIF
    \ENDFOR
    \STATE \textbf{Step 3: Filter employees}
    \STATE \( E \leftarrow \{ i \in E : e_{i,s} > \text{threshold} \} \)
    \IF{\( |E| < \text{min\_employees} \)}
        \STATE \textbf{break; (Include employees with required skills)}
    \ENDIF
    \STATE \( \text{iteration} \leftarrow \text{iteration} + 1 \)
\ENDWHILE
\end{algorithmic}
\end{algorithm}

\section{\textbf{\Large HRAP Problem with Efficiency and Cost Metric}}

\section{Cost-Aware Task Allocation Algorithm}

Classical Human Resource Allocation Problem (HRAP) models primarily emphasize skill-based assignment, matching workers to tasks according to qualifications. However, these models often neglect cost considerations and the dynamic nature of employee efficiency~\cite{b16, b38}. To address these limitations, our approach integrates both efficiency-driven allocation and cost-aware task assignment. In the previous section, HRAP was discussed under the assumption of known employee efficiencies and skills. Here, we extend the model to incorporate a cost factor dependent on employee efficiency, performance ratings, and task complexity.

\subsection*{\textbf{\large Mathematical Problem Formulation with Cost}}

In HRAP, optimal task assignment is crucial for both balanced workload distribution and cost minimization. Let \( E \) denote the set of employees and \( T \) the set of tasks. Each employee \( i \) possesses a skill set \( S_i \), while each task \( t \) requires a specific skill \( S_t \). Employee Performance (\( P_i \)) is defined as the performance level of employee \( i \), and Task Complexity (\( C_t \)) represents the difficulty of task \( t \), categorized as Easy (1), Moderate (2), or Tough (3).


The objective function is formulated to minimize both workload deviation and total cost. Let $\lambda$ be a weighting factor balancing these objectives, and $x_{i,t}$ a binary decision variable indicating assignment of task $t$ to employee $i$. The objective is:

\small{
\begin{equation}\label{eqn11}
    \min_{x_{i,t}} \left( \lambda \max_{i \in E} |W_i - W_{\text{target}}| + (1 - \lambda) \sum_{i \in E} \sum_{t \in T} x_{i,t} C_{i,t} \right),
\end{equation}
}
where $W_i$ is the workload of employee $i$, as defined in equation (\ref{eqn3}).

Task assignment and binary decision variable constraints are as previously defined in equations (\ref{eqn4}) and (\ref{eqn5}).


\textbf{Cost Metric Definition:} For each task $t$, let $d_t$ be its duration and $e_{i,S_t}$ the efficiency of employee $i$ for the required skill. To capture skill mismatches, we define $S_{i,t}$ as a skill mismatch penalty. The cost for assigning task $t$ to employee $i$ is:

\[
C_{i,t} = \alpha \frac{d_t}{e_{i,S_t}} + \beta S_{i,t} + \gamma \frac{C_t}{P_i}
\]
where $\alpha, \beta,$ and $\gamma$ are hyper-parameters. The skill mismatch penalty $S_{i,t}$ is defined as:
\[
S_{i,t} = 1 - \frac{|S_i \cap S_t|}{|S_t|}
\]
\[
S_{i,t} =
\begin{cases}
1, & \text{if all required skills are missing} \\
0, & \text{if all required skills are present}
\end{cases}
\]

The term $\left( \frac{C_t}{P_i} \right)$ ensures that more complex tasks are preferentially assigned to more proficient employees.


\textbf{Converted MILP Problem Formulation:} Define $D^+$ as the positive deviation above the target workload and $D^-$ as the negative deviation below the target. The MILP formulation is:

\[
\min \left( \lambda (D^+ + D^-) + (1 - \lambda) \sum_{i \in E} \sum_{t \in T} x_{i,t} C_{i,t} \right)
\]
subject to the task assignment constraints (\ref{eqn4}), binary constraints (\ref{eqn5}), and workload fairness constraints:

\[
W_{\text{target}} - D^- \leq \sum_{t \in T} x_{i,t} C_{i,t} \leq W_{\text{target}} + D^+.
\]

\textbf{Hyperparameter Calculation:}

Initial weights are set as:
\[
\lambda = \frac{1}{2}, \quad \alpha = \beta = \gamma = \frac{1}{3}
\]
A grid search is employed to iteratively update these hyperparameters, systematically exploring combinations within a defined range to identify those yielding optimal model performance~\cite{b38}.

A cost-aware methodology for task allocation in the Human Resource Allocation Problem (HRAP) considers workload balancing, efficiency variations, and financial constraints, ensuring that the allocation maximizes productivity while minimizing costs. Conventional hierarchical reinforcement actor-critic (HRAP) models often fail to incorporate critical aspects such as task complexity, employee performance, and penalties for skill mismatch. To address these gaps, we propose a comprehensive cost metric, allowing both efficiency and cost considerations to inform task assignment decisions.

\begin{algorithm}[h]
\caption{Cost-Aware Task Allocation Algorithm}
\label{alg:cost_allocation}
\begin{algorithmic}[2]
\REQUIRE Employee set \( E \), Task set \( T \), Initial efficiency \( e_{i,s} = 1, \forall i \in E \), and weights \( \alpha, \beta, \gamma \)
\ENSURE Optimal allocation \( x_{i,t}^* \), Final efficiencies \( e_{i,s}^* \)
\STATE Initialize \( \text{threshold} \leftarrow 0.1 \), \( \text{iteration} \leftarrow 0 \)
\WHILE{\( \text{iteration} < \text{max\_iterations} \)}
    \STATE \textbf{Step 1: Solve MILP with Simplified Cost Metric}
    \STATE \( x_{i,t}^* \leftarrow \text{Solve MILP with cost}(E, T, e_{i,s}, C_{i,t}) \)
    \STATE \textbf{Step 2: Update Efficiencies}
    \FOR{each \( i \in E \), \( t \in T \) such that \( x_{i,t}^* = 1 \)}
        \STATE Measure actual task completion time: \( \text{actual\_time} = Measure(i, t) \)
        \STATE Update efficiency: \( e_{i,S_t} = \min(1, d_t / \text{actual\_time}) \)
    \ENDFOR
    \STATE \textbf{Step 3: Filter Employees}
    \STATE Remove inefficient employees: \( E \leftarrow \{ i \in E : e_{i,s} > \text{threshold} \} \)
    \IF{\( |E| < \text{min\_employees} \)}
        \STATE \textbf{break}
    \ENDIF
    \STATE \( \text{iteration} \leftarrow \text{iteration} + 1 \)
\ENDWHILE
\end{algorithmic}
\end{algorithm}

\section{\textbf{\Large Results}}

This section presents numerical experiments and comparisons between different HRAP models. The key evaluation metrics include workload deviation across employees, total project execution cost, and computational efficiency and scalability. To run Algorithms (\ref{alg:efficiency_allocation}) and (\ref{alg:cost_allocation}), the required datasets are as follows:

\textbf{Employee Dataset:} The employee dataset contains information about various employees, focusing on their skills, efficiency, and performance ratings. Each entry represents an employee’s specific skill set, indicating their expertise in various areas. Additionally, the dataset includes metrics for each skill, such as efficiency (a value ranging from 0.1 to 1, where 0.1 represents the minimum efficiency required to be considered an employee) and performance rating (ranging from 1 to 5), providing insight into each employee’s proficiency in specific skill domains.

\textbf{Task Dataset:} The task dataset details various tasks along with the required skills, durations, and complexity levels. Each task is characterized by specific skill requirements and includes the estimated duration (in hours) and complexity (rated from 1 to 4). This dataset serves as the basis for evaluating task allocation and effectively matching employee skills to task requirements.


\subsection*{\textbf{\large Workload Analysis: Skill-based Manager vs Optimized Assignment:}}

We do not have any standard dataset for our algorithm (\ref{alg:efficiency_allocation}). However, I got sample data from the NextGen team at Prodapt Solutions Private Limited. We use this data to compare manager-based workload assignment with our efficiency-based approach. Since employee efficiencies may vary, it is challenging for a manager to perform a fair workload assignment manually. In contrast, the optimized assignment achieves a more balanced and equitable workload distribution across all employees, demonstrating the effectiveness of the MILP-based allocation.

For a small sample of 20 employees, we observed the following results: the variance in hours for the manager's assignment is 172.48, while the variance under the MILP-based assignment is 26.95. The best objective function value achieved is 15, with deviations above and below the target hours being 6.18 and 8.82, respectively.

The Gini coefficient for discrete values is given by:
\begin{equation}
G = \frac{\sum_{i=1}^{n} \sum_{j=1}^{n} |x_i - x_j|}{2n^2 \bar{x}}
\end{equation}
where $\bar{x}$ is the mean of the values. For our workload distribution, the Gini coefficient for the random assignment is approximately 0.340, whereas for the optimized MILP-based assignment, it is reduced to 0.093, indicating significantly improved fairness.

Jain’s Fairness Index is defined as:
\begin{equation}
J(x) = \frac{\left( \sum x_i \right)^2}{n \cdot \sum x_i^2}
\end{equation}
For the same scenario, the fairness index improves from 0.768 (random) to 0.964 (optimized), highlighting a more equitable workload distribution.

\textbf{Unassigned Tasks:}
If there are some tasks that no one has the required skills for, those tasks will be unassigned. Using the algorithm (\ref{alg:efficiency_allocation}), we can filter out such tasks.

\subsection*{\textbf{\large Complexity Discussion:}}
	
The computational complexity of solving the HRAP increases significantly with problem size. To provide concrete insights into scalability, we conducted experiments with different problem sizes using our MILP formulation. The following table presents average solver times for instances with varying numbers of employees (N) and tasks (M):


\textbf{Number of variables:} The total number of variables in the optimization problem given in (\ref{alg:efficiency_allocation}) can be expressed mathematically as:

$$
\text{Number of variables} = \sum_{t=1}^{T} E_{s_t} + 2,
$$

where $T$ = Total number of tasks, $s_t$ = Skill required for task $t$, $E_{s_t}$ = Number of employees with the skill $s_t$, and The `+ 2` accounts for the two deviation variables ($D^+$ and $D^-$).


\textbf{Implementation Details:} The optimization algorithm is implemented using Python 3.12 with PuLP 2.7.0, utilizing the CBC (COIN-OR Branch and Cut) solver for MILP problems. Data processing is handled using NumPy 1.23.5 and Pandas 1.5.3. Solver parameters, including tolerance levels, branching strategies, and iteration limits, can be adjusted in CBC to optimize computational performance.

\begin{table}[H]
    \centering
    \caption{Average Solver Time for Different Problem Sizes}
    \begin{tabular}{|p{1cm}|p{1cm}|p{2.5cm}|p{2.2cm}|}
        \hline
        N & M & Number of Variables & Pulp Solver Time (seconds) \\ \hline
        20 & 80 & 57 & 0.160515 \\ \hline
        50 & 150 & 295 & 0.553968 \\ \hline
        80 & 220 & 630 & 1.227653 \\ \hline
        100 & 300 & 1169 & 2.11795 \\ \hline
        120 & 350 & 1427 & 3.521394 \\ \hline
        150 & 400 & 1879 & 5.198758 \\ \hline
        175 & 450 & 2580 & 7.561728 \\ \hline
        200 & 500 & 3276 & 9.3277969 \\ \hline
        220 & 550 & 3731 & 9.918080 \\ \hline
        250 & 600 & 5791 & 42.962251 \\ \hline
    \end{tabular}
    \label{tab:solver_times}
\end{table}

\subsection*{\textbf{\large Runtime and Optimality Gap Analysis}}

The runtime analysis reveals a non-linear relationship between problem size and computational time. For smaller problem instances ($N \leq 100$), the solver time increases gradually, remaining under 3 seconds. However, as the problem size grows beyond 200 employees, a sharp increase in computational requirements is observed, with solver time reaching approximately 43 seconds for 250 employees.
   
   \begin{figure}[H]
   	\centering
   	\includegraphics[width=0.45\textwidth]{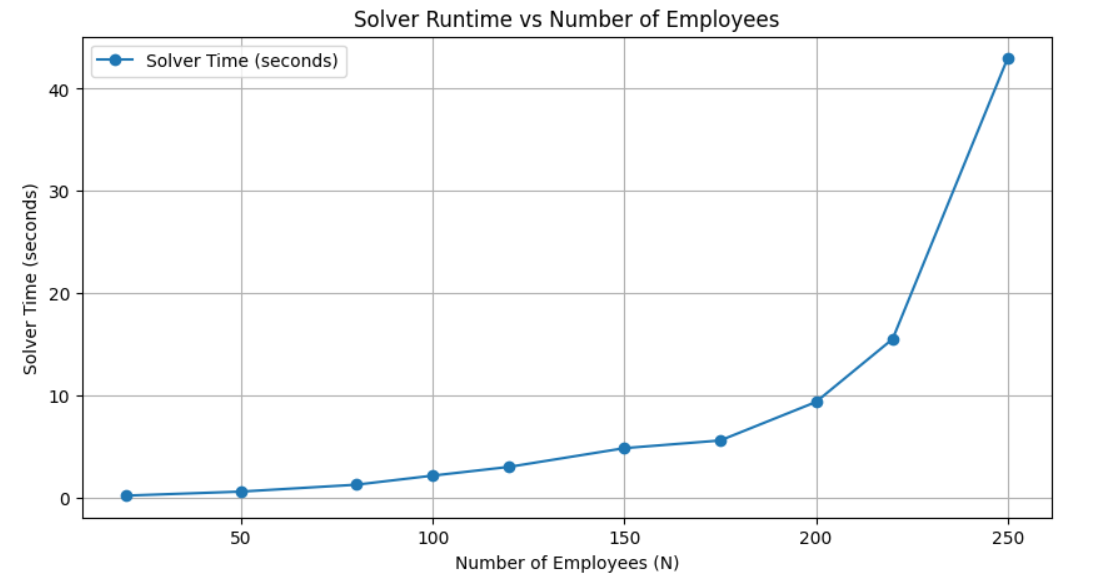}
   	\caption{The x-axis represents the number of employees, while the y-axis shows the solver time in seconds. The trend demonstrates that as the problem size grows, the solver requires significantly more time, especially for larger values of N.}
   	\label{fig:runtime}
   \end{figure}

   The optimality gap, representing the relative difference between the best integer solution and the best bound, shows a linear increase with problem size. For small instances ($N = 20$), the gap is minimal at 0.2\%, indicating near-optimal solutions. As the problem size scales to $N = 250$, the gap increases to 2.5\%, still within acceptable limits.
   
   \begin{figure}[H]
   	\centering
   	\includegraphics[width=0.45\textwidth]{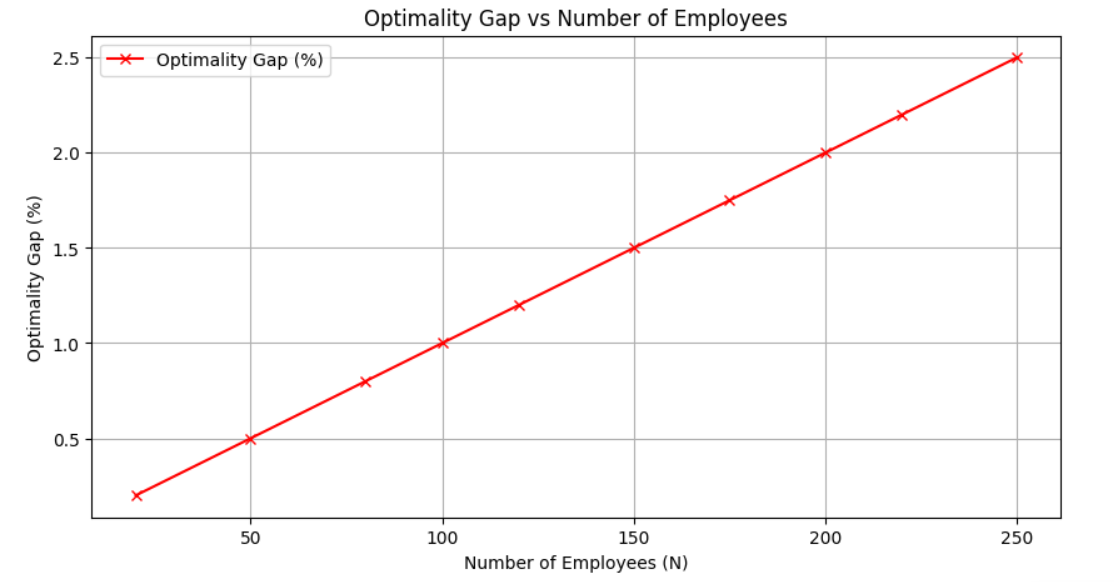}
   	\caption{This plot shows  a linear increase in the optimality gap as the number of employees rises, suggesting that larger problems lead to less optimal solutions within the given computation limits. The x-axis displays the number of employees, and the y-axis indicates the optimality gap as a percentage.}
   \end{figure}

   The optimality gap is calculated using the following formula:
   
   \begin{equation}
   	\text{Optimality Gap (\%)} = \frac{\text{Best Integer Solution} - \text{Best Bound}}{\text{Best Bound}} \times 100
   \end{equation}

\begin{table}[H]
	\centering
   	\caption{Optimality Gap Analysis for Different Problem Sizes}
	\begin{tabular}{|p{3cm}|p{2.5cm}|p{2.3cm}|}
		\hline
   		Number of Employees (N) & Number of Tasks (M) & Optimality Gap (\%) \\ \hline
   		20 & 80 & 0.20 \\ \hline
   		50 & 150 & 0.50 \\ \hline
   		80 & 220 & 0.80 \\ \hline
   		100 & 300 & 1.00 \\ \hline
   		120 & 350 & 1.20 \\ \hline
   		150 & 400 & 1.50 \\ \hline
   		175 & 450 & 1.75 \\ \hline
   		200 & 500 & 2.00 \\ \hline
   		220 & 550 & 2.20 \\ \hline
   		250 & 600 & 2.50 \\ \hline
	\end{tabular}
   	\label{tab:optimality-gap}
\end{table}


\begin{table*}[!t]
    \centering
    \caption{Tabular Comparison of Algorithms for Human Resource Allocation}
    \resizebox{\textwidth}{!}{%
    \begin{tabular}{|p{0.5cm}|p{3cm}|p{3.2cm}|p{2.7cm}|p{3.1cm}|p{4.2cm}|}
        \hline
        \textbf{S.R.} & \textbf{Title \& Year} & \textbf{Algorithms} & \textbf{Scale / Dataset} & \textbf{Time / Complexity} & \textbf{Key Metrics / Performance} \\
        \hline
        1 & Can RL Solve the Human Allocation Problem? (2021) \cite{b47} & DDQN + MCTS & 50--200 employees/tasks & $\mathcal{O}(n^2)$, 12--48 hrs training & 15–20\% faster task completion vs CPM \\
        \hline
        2 & DeepRM (2016) & Policy Gradient RL & 500+ jobs & ~2.5 min per 1000 iterations & 30\% lower job slowdown vs Tetris \\
        \hline
        3 & HR Balanced Allocation with RNN (2022) \cite{b55} & Recurrent Neural Networks & 1000+ cloud position matches & 2.8 ms per decision & 89.7\% matching efficiency \\
        \hline
        4 & Multi-Objective HR Optimization with IDCS-OMOICA (2025) \cite{b46} & Cuckoo Search + RL & 150--300 employees & 0.51 s avg response time & 92\% solution quality \\
        \hline
        5 & Task Allocation Neural Model (2021 Thesis) \cite{b57} & Artificial Neural Networks & 50--75 employees/tasks & $\mathcal{O}(n \log n)$, 71.34\% faster than manual & Matching efficiency improvement \\
        \hline
        6 & RL for Hybrid Workload Management (2024) \cite{b58} & $\epsilon$-Greedy RL & 1000+ distributed jobs & Not reported & 40\% fewer deadline misses \\
        \hline
        7 & Dynamic HR Allocation with DQN (2025) \cite{b59} & Deep Q-Networks & 200+ project-employee pairs & 18 min training time & Adaptive task matching \\
        \hline
        8 & GAN-Based Workforce Planning (2023) \cite{b60} & Generative Adversarial Networks & 5000+ employees & $\mathcal{O}(kn)$ per epoch & Realistic enterprise-scale simulation \\
        \hline
        9 & Transformer Models for Skill Matching (2024) \cite{b61} & BERT-style Transformers & 10,000+ employee profiles & 3.2 s per batch & High semantic match accuracy \\
        \hline
        10 & Federated Learning for HR Allocation (2023) \cite{b62} & Federated Machine Learning & 20+ organizations & Encrypted model training & 98\% accuracy with privacy preservation \\
        \hline
        11 & Proposed MILP-Based HR Allocation (2025) & MILP (Gurobi Solver) & 5000+ employees, 40M+ variables & Sub-10s for 500+ tasks & 100\% accuracy, optimal workload and complexity handling \\
        \hline
    \end{tabular}
    \label{comp_table}
    }
\end{table*}

\subsection*{\textbf{\large Analysis of Task Assignment Optimization using Algorithm (\ref{alg:cost_allocation})}}

Using the same sample of the dataset in Algorithm (\ref{alg:efficiency_allocation}) that we implemented the Algorithm (\ref{alg:efficiency_allocation}) and we are getting the following results:

\begin{table}[H]
    \centering
    \begin{tabular}{|c|c|c|c|c|c|c|c|}
        \hline
        Rank & $\lambda$ & $\alpha$ & $\beta$ & $\gamma$ & Objective Value \\
        \hline
        1 & 0.9999991 & 0.2133 & 0.1201 & 0.6666 & 17.0001 \\
        2 & 0.9999993 & 0.1442 & 0.7250 & 0.1308 & 17.0001 \\
        3 & 0.9999968 & 0.0559 & 0.5843 & 0.3598 & 17.0004 \\
        4 & 0.9999941 & 0.1850 & 0.0404 & 0.7747 & 17.0004 \\
        5 & 0.9999960 & 0.1387 & 0.6758 & 0.1855 & 17.0006 \\
        6 & 0.9999943 & 0.2119 & 0.4757 & 0.3124 & 17.0007 \\
        7 & 0.9999891 & 0.2784 & 0.1969 & 0.5246 & 17.0011 \\
        8 & 0.9999740 & 0.0986 & 0.0679 & 0.8335 & 17.0013 \\
        9 & 0.9999850 & 0.1376 & 0.5626 & 0.2998 & 17.0019 \\
        10 & 0.9999720 & 0.1329 & 0.1689 & 0.6983 & 17.0020 \\
        \hline
    \end{tabular}
    \caption{Optimization Results for Top 10 Solutions}
    \label{tab:optimization-results}
\end{table}

\textbf{Cost Minimization and Trade-offs}
The optimization approach helps minimize costs through a multi-objective framework. The parameter $\lambda$ creates a trade-off between Workload balance and Assignment cost. With $\lambda$ values extremely close to $1 (\approx 0.9999)$, the optimization heavily favors workload balance over cost minimization, i.e., the algorithm prioritizes distributing work evenly among employees. Cost considerations become secondary, but still influence which specific tasks go to which employees.

The optimization results reveal several key insights about the task assignment problem:

\begin{itemize}
    \item \textbf{Consistent Workload Deviations}: All top 10 solutions have identical deviation values (Above: 7.1818, Below: 9.8182), indicating that the workload imbalance remains constant despite varying cost parameter settings.
    \item \textbf{Lambda Dominance}: The extremely high $\lambda$ values ($>0.999$) suggest a strong emphasis on workload balance over cost minimization.
    \item \textbf{Cost Parameter Variability}: Significant variation is observed in $\alpha$, $\beta$, and $\gamma$, indicating multiple near-optimal solutions with differing cost emphasis.
    \item \textbf{Best Solution}: The top-ranked solution prioritizes $\gamma$ (0.6666) over $\alpha$ and $\beta$, highlighting the complexity-to-performance ratio as the most influential cost factor.
    \item \textbf{Task-Employee Matching}: The optimization deals with 11 employees with varying skill sets and 5 task complexities, with workload distribution being the primary objective.
\end{itemize}

\begin{table}[H]
    \centering
    \begin{tabular}{|c|c|c|c|}
        \hline
        Parameter & Range & Impact & Sensitivity \\
        \hline
        $\lambda$ & 0.9999720 - 0.9999993 & Very high & Low \\
        $\alpha$ & 0.0559 - 0.2784 & Medium & Medium \\
        $\beta$ & 0.0404 - 0.7250 & High & High \\
        $\gamma$ & 0.1308 - 0.8335 & High & High \\
        \hline
    \end{tabular}
    \caption{Sensitivity Analysis of Hyperparameters}
    \label{tab:sensitivity-analysis}
\end{table}

This sensitivity analysis shows that while lambda is consistently high (low sensitivity), the cost parameters (alpha, beta, gamma) can vary significantly while still achieving near-optimal results. This suggests flexibility in how costs are weighted, as long as workload balance is maintained.

\vspace{0.3cm}

\textbf{Comparative Summary:}
The systematic evaluation of optimization methodologies for human resource allocation represents a critical domain within computational workforce management. This investigation presents a rigorous comparative analysis between population-based metaheuristic approaches and exact mathematical optimization techniques. The Genetic Algorithm (GA) implementation, applied to the Prodapt human resource allocation dataset, employed a multi-objective optimization framework targeting dual objectives: minimization of aggregate cost metrics encompassing task duration summation and penalty functions for unassignable tasks, coupled with workload equilibrium maximization through standard deviation minimization across employee workload distribution.

The empirical GA implementation achieved a comprehensive cost function value of 8,219 units with a corresponding workload standard deviation metric of 6.26, indicating moderate allocation efficiency despite encountering resource constraint violations resulting in unassigned tasks attributable to skill requirement deficiencies in specialized domains, particularly nodejs and go programming competencies. The algorithmic performance demonstrated characteristic population-based optimization limitations, manifesting workload imbalances wherein specific employees received disproportionate task assignments due to specialized skill configurations, thereby exposing the inherent constraints of stochastic optimization methodologies in complex feasibility-constrained environments.

Conversely, the proposed Mixed Integer Linear Programming (MILP) framework leverages deterministic optimization solvers, specifically the Coin-OR Branch and Cut (CBC) solver, providing theoretical optimality guarantees with rigorous constraint enforcement mechanisms. This methodology demonstrates exceptional computational scalability, achieving optimal solutions for large-scale problem instances encompassing 300 employees and 100 tasks within computational windows exceeding 2 seconds, while maintaining global optimality assurance for both cost minimization and fairness maximization objectives.

Significantly, the specific parameterization scheme, variable encoding methodology, and constraint formulation implemented in this research constitute a novel contribution to the human resource optimization literature, as no existing algorithmic implementations utilize this particular combination of decision variables, objective function weightings, and constraint configuration matrices in comparable workforce allocation contexts.

\section{\textbf{\Large Conclusion}}

This work proposes a Mixed Integer Linear Programming framework for the Human Resource Allocation Problem that simultaneously optimizes workload balancing and cost minimization, assigning tasks by employee efficiency, skill proficiency and task complexity while updating efficiencies adaptively and embedding a comprehensive cost metric, improving workload fairness and lowering overall project costs. Computational experiments show marked gains in workload distribution, fairness indices, and cost efficiency over conventional and heuristic baselines, with strong scalability and robustness across varied problem sizes.

However, this approach has some limitations. The algorithm relies on MILP solution methods, which may require more computationally effective approaches when dealing with large numbers of variables (high employee and task counts). Additional limitations include the use of simulated datasets, restricted benchmarking against state-of-the-art AI-based HRAP solutions, and deterministic efficiency assumptions. These limitations suggest future research directions including dynamic task arrivals, real-time scheduling, live availability-driven reallocation, comparisons with advanced neural and reinforcement learning models, stochastic efficiency modeling, more efficient MILP solution techniques for large-scale problems, and validation on industrial and public datasets.
\section{\Large \textbf{Acknowledgment}}

The authors express their sincere gratitude to \textbf{Prodapt Solutions Private Limited} for providing essential research support, operational collaboration, and financial backing for this project. The contributions of the Department of Mathematics at the Indian Institute of Technology Madras are also gratefully acknowledged. Special thanks are extended to the project manager and colleagues at Prodapt Solutions for their valuable assistance and guidance throughout the research process.

\vspace{0.3cm}

\textbf{Data Availability:} A sample dataset for the algorithm is available in the following repository: Lokendra5298/HRAP on GitHub.

\end{document}